 \documentclass[final,authoryear,3p,times]{elsarticle}

 \usepackage{graphicx}
 \usepackage{subfig} 
\usepackage{color}
\usepackage{amssymb, amsmath}
\usepackage{amsthm}
\usepackage{bbm}

\journal{Statistics \& Probability Letters}

\newcommand{\ind}[1]{\mathbbm{1}_{\{#1\}}}
\newcommand{\indnex}{\mathbbm{1}_{\overline{\cE}}}
\newcommand{\bs}[1]{\mathbf{#1}}

\newcommand{\dE}{\mathbb{E}}
\newcommand{\dF}{\mathbb{F}}
\newcommand{\dG}{\mathbb{G}}
\newcommand{\dP}{\mathbb{P}}
\newcommand{\dR}{\mathbb{R}}
\newcommand{\dT}{\mathbb{T}}

\newcommand{\wh}{\widehat}

\newcommand{\veps}{\varepsilon}

\newcommand{\cE}{\mathcal{E}}
\newcommand{\cF}{\mathcal{F}}
\newcommand{\cG}{\mathcal{G}}

\newcommand{\cL}{\mathcal{L}}
\newcommand{\cN}{\mathcal{N}}

\newcommand{\reff}[1]{(\ref{#1})}

\def\build#1_#2^#3{\mathrel{\mathop{\kern 0pt#1}\limits_{#2}^{#3}}}
\def\liml{\build{\longrightarrow}_{}^{\cL}}

\begin{document}

\newtheorem{lemma}{Lemma}[section]
\newtheorem{proposition}[lemma]{Proposition}
\newtheorem{theorem}[lemma]{Theorem}
\newtheorem{corollary}[lemma]{Corollary}

\begin{frontmatter}



\title{Asymmetry tests for Bifurcating Auto-Regressive Processes  with missing data}


\author[ad1]{Beno\^\i te de~Saporta}
\author[ad2]{Anne G\'egout-Petit}
\author[ad3]{Laurence Marsalle}
\address[ad1]{Univ. Bordeaux, GREThA CNRS UMR~5113, IMB CNRS UMR~5251, \\
CNRS, GREThA UMR~5113, IMB UMR~5251,\\
INRIA~Bordeaux~Sud~Ouest team CQFD,\\ 
351 cours de la Lib\'eration, 33405 Talence Cedex, France,\\ saporta@math.u-bordeaux1.fr,\\
 tel: +33 5 24 57 41 69, fax:  +33 5 24 57 40 24}
\address[ad2]{Univ. Bordeaux, IMB  CNRS UMR~5251,\\
CNRS, IMB UMR~5251,\\
INRIA~Bordeaux~Sud~Ouest team CQFD,\\
351 cours de la Lib\'eration, 33405 Talence Cedex, France,\\
anne.petit@u-bordeaux2.fr}
\address[ad3]{Univ. Lille 1, Laboratoire Paul Painlev\'e, CNRS UMR~8524,\\
CNRS, Laboratoire Paul Painlev\'e, UMR~8524
59655 Villeneuve d'Ascq Cedex, France,\\
laurence.marsalle@univ-lille1.fr}

\begin{abstract}
We present symmetry tests for bifurcating autoregressive processes (BAR) when some data are missing. BAR processes typically model cell division data. Each cell can be of one of two types \emph{odd} or \emph{even}. The goal of this paper is to study the possible asymmetry between odd and even cells in a single observed lineage. We first derive asymmetry tests for the lineage itself, modeled by a two-type Galton-Watson process, and then derive tests for the observed BAR process. We present applications on both simulated and real data.
\end{abstract}

\end{frontmatter}



\section{Introduction}
\label{s:intro}
Bifurcating autoregressive processes (BAR) were first introduced by \citet{CoSt86}. They generalize autoregressive processes when data are structured as a binary tree. Typically, they are involved in statistical studies of cell lineages, since each cell in one generation gives birth to two offspring in the next one. {These two offspring sometimes need to be distinguished according to some biological property, leading to the notion of \emph{type}.} Cell lineage data consist of observations of some quantitative characteristic of the cells (e.g. their growth rate) over several generations descended from an initial cell. 
More precisely, the initial cell is labelled $1$, and the two offspring of cell $k$ are labelled $2k$ and $2k+1${, where $2k$ stands for one type, thus called \emph{even}, and $2k+1$ for the other type, thus called \emph{odd}}. If $X_k$ denotes the quantitative characteristic of cell $k$, then the first-order asymmetric BAR process is given by
\begin{equation}\label{defbar}
\left\{
    \begin{array}{lcccccl}
     X_{2k} & = & a &+ &bX_k &+ &\varepsilon_{2k}, \\
     X_{2k+1} & = & c & +& dX_k &+ &\varepsilon_{2k+1},
    \end{array}\right.
\end{equation}
for all $k\geq 1$. The noise sequence $(\varepsilon_{2k}, \varepsilon_{2k+1})$ represents environmental effects, while $a,b,c,d$ are unknown real parameters related to the inherited effects. They are allowed to be different for the odd and even sisters.

\smallskip

Various estimators are studied in the literature for the unknown parameters $a,b,c,d$, see
\citep{BSG09, DM08, Guy07}. This paper derives further properties of the estimators given in \citep{SGM11}, where the genealogy is modeled by a two-type Galton Watson process (GW), allowing the reproduction laws to depend on both the mother's and daughter's types. Indeed, the aim of this paper is to propose asymmetry tests for both the two-type GW process defining the genealogy of the cells, as well as the BAR process with missing data. More precisely, we first study the difference of the means of the reproduction laws for even and odd mother cells in the GW process. Then we investigate the difference of the parameters $a,b$ and $c,d$ and the difference between the fixed points $a/(1-b)$ and $c/(1-d)$. 
We propose Wald's type tests based on the asymptotic normality of the various estimators. A detailed study on simulated data, as well as a new investigation of the \emph{Escherichia coli} data of \citet{SMT05} are provided.

\smallskip

The paper is organized as follows. We start with introducing in Section~\ref{s:notation} some notation that will be used throughout the paper. In Section~\ref{s:lineage}, we derive Wald's tests for the two-type GW process describing the genealogy of the cells. In Section~\ref{s:BAR}, we derive Wald's test for the BAR process, to test asymmetry of its parameters. In section~\ref{s:simu} we apply our tests to simulated data. Finally, in Section~\ref{s:Ecoli} we apply our tests to \citep{SMT05} data.
\section{Notation}
\label{s:notation}
For all $n\geq 1$, denote the n-th generation by
$\dG_n = \{k \in \mathbb{N} \vert  2^n\leq k \leq 2^{n+1}-1\}$.
Denote by $\dT_n = \bigcup_{\ell=0}^n\dG_{\ell}$,
the sub-tree of all cells 
up to the $n$-th generation. The cardinality $|\dG_n|$ of $\dG_n$  is $2^n$, while that of $\dT_n$ is $|\dT_n|=2^{n+1}-1$. 
We need to distinguish the cells in $\dG_n$ and $\dT_n$ according to their \emph{type}. The type \emph{even} will be labelled \emph{type 0} and the type \emph{odd} will be labelled \emph{type 1}. We set 
$ \dG_n^0=\dG_n \cap (2 \mathbb{N})$, $\dG_n^1=\dG_n \cap (2 \mathbb{N} + 1)$, $\dT_n^0=\dT_n \cap (2 \mathbb{N})$, $\dT_n^1=\dT_n \cap (2 \mathbb{N}+1)$.
We encode the presence or absence of the cells by the process $(\delta_k)$: if cell $k$ is observed, then $\delta_k=1$, if cell k is missing, $\delta_k=0$. We define the sets of observed cells as 
$\dG_n^*=\{k \in \dG_n : \delta_k =1\}$ and $\dT_n^*=\{k \in \dT_n : \delta_k =1\}$.
Finally, let $\mathcal{E}$ be the event corresponding to the cases when there are no cell left to observe in the current generation: $\mathcal{E} = \bigcup_{n \ge 1} \{|\dG_n^*|=0 \}$ and $\overline{\cE}$ the complementary set of $\cE$. For $n \ge 1$, we define the number of observed cells among the n-th generation, distinguishing according to their type: $Z_n^0=|\dG_n^* \cap 2 \mathbb{N}|$, $Z_n^1=|\dG_n^* \cap (2 \mathbb{N}+1)|$,
and we set, for all $n \ge 1$, $\bs{Z}_n=(Z_n^0,Z_n^1)$. Note that for $i\in\{0,1\}$ and $n\geq 1$ one has
$Z_n^i=\sum_{k\in\dG_{n-1}}\delta_{2k+i}$.
\section{Asymmetry in the lineage}
\label{s:lineage}
We now describe the mechanism generating the observation process $(\delta_k)$. The resulting process  $(\bs{Z}_n)$  is a two-type GW process. We recall some assumptions similar to \citep{SGM11} and mostly taken from \citep{Har63}.
\subsection{Model and assumptions}
\label{ss:GWmodel}
We define the cells genealogy by a two-type GW process $(\bs{Z}_n)$. All cells reproduce independently and with a reproduction law depending only on their type. For a mother cell of type $i\in\{0,1\}$, we denote by $p^{(i)}(j_0,j_1)$ the probability that it has $j_0$ daughter of type $0$ and $j_1$ daughter of type $1$. 
For the cell division process we are interested in, one clearly has 
$p^{(i)}(0,0)+p^{(i)}(1,0)+p^{(i)}(0,1)+p^{(i)}(1,1)=1$.
The reproduction laws have moments of all order, and we can thus define the descendants matrix $\bs{P}$
\begin{equation*}
\bs{P}=\left(\begin{array}{cc}
           p_{00} & p_{01} \\
           p_{10} & p_{11}
          \end{array}\right),
\end{equation*}
where $p_{i0} = p^{(i)}(1,0) + p^{(i)}(1,1)$ and $p_{i1} = p^{(i)}(0,1) + p^{(i)}(1,1)$, for $i \in \{0,1\}$: $p_{ij}$ is the expected number of descendants of type $j$ of a cell of type $i$. It is well-known that when all the entries of the matrix $\bs{P}$ are positive, $\bs{P}$ has a strictly dominant positive eigenvalue, denoted $\pi$, which is also simple. We make the following main assumption.
\begin{description}
 \item[(\bf{AO})] All entries of the matrix $\bs{P}$ are positive and the dominant eigenvalue is greater than one: $\pi > 1$ .
\end{description}
In this case, there exist left and right eigenvectors for $\pi$ which are component-wise positive. Let $\bs{z}=(z^0,z^1)$ be such a left eigenvector satisfying $z^0+z^1=1$.{This dominant eigenvalue $\pi$ is related to the extinction of the process: assumption (\textbf{AO}) means that the GW process $(\bs{Z}_n)$ is super-critical, and ensures that extinction is not almost sure: $\dP(\mathcal{E}) < 1$. Besides, on $\overline{\cE}$, $|\dT_n^*|^{-1}\sum_{l=1}^nZ_l^i$ converges to $z^i$, meaning that $z^i$ is the asymptotic proportion of cells of type $i$ in a given sub-tree.}
\subsection{Asymmetry test}
\label{ss:GWtest}
We first propose estimators for the parameters of the GW process and study their asymptotic properties. Our context is very specific because the available information given by $(\delta_k)$ is more precise than the one given by $(\bs{Z}_n)$ usually used in the literature, see e.g. \citep{Gut91}. The empiric estimators of the reproduction probabilities using data up to the $n$-th generation are then, for $i,j_0,j_1$ in $\{0,1\}$
\begin{equation*}
\widehat{p}_n^{(i)}(j_0,j_1) = \frac{\sum_{k\in\dT_{n-2}}\delta_{2k+i}\phi_{j_0}(\delta_{2(2k+i)})\phi_{j_1}(\delta_{2(2k+i)+1})}{\sum_{k\in\dT_{n-2}}\delta_{2k+i}},
\end{equation*}
where $\phi_0(x)=1-x$, $\phi_1(x)=x$, and if the denominators are non zero (the estimation is zero otherwise). The strong consistency is easily obtained on the non-extinction set $\overline{\cE}$ e.g. by martingales methods as in \citep{SGM11}.
\begin{lemma}\label{lem:CVp}
Under \emph{(\textbf{AO})} and for all $i$, $j_0$ and $j_1$ in $\{0,1\}$, one has
$\lim_{n\rightarrow\infty}\ind{|\dG_n^*|>0}\widehat{p}_n^{(i)}(j_0,j_1)=\mathbbm{1}_{\overline{\cE}}p^{(i)}(j_0,j_1)$ a.s.
\end{lemma}
Set $\bs{p}^{(i)}=(p^{(i)}(0,0),p^{(i)}(1,0),p^{(i)}(0,1),p^{(i)}(1,1))^t$ the vector of the $4$ reproduction probabilities for a mother of type $i$, $\bs{p}=((\bs{p}^{(0)})^t, (\bs{p}^{(1)})^t)^t$ the vector of all $8$ reproductions probabilities and $\widehat{\bs{p}}_{n}=(\wh{p}_n^{(0)}(0,0),\ldots,\wh{p}_n^{(1)}(1,1))^t$ its estimator. As $\dP(\overline{\cE})\neq 0$, we define the conditional probability $ \dP_{\overline{\cE}}$ by
$\dP_{\overline{\cE}}(A) = {\dP(A \cap \overline{\cE})}/{ \dP(\overline{\cE})}$ for all event $A$.
\begin{theorem}\label{th:TCLGW}
Under assumption \emph{(\textbf{AO})}, we have the convergence 
$\sqrt{|\dT_{n-1}^*|}(\widehat{\bs{p}}_{n}-\bs{p})
\liml
\cN(0, \bs{V} )$ on $(\overline{\cE},\mathbb{P}_{\overline{\cE}})$,
with
$$\bs{V}=\left(
\begin{array}{cc}
\bs{V}^0/z^0&0\\
0&\bs{V}^1/{z^1}
\end{array}\right)$$
and for all $i$ in $\{0,1\}$, $\bs{V}^i=\bs{W}^i-\bs{p}^{(i)}(\bs{p}^{(i)})^t$, $\bs{W}^i$ is a $4\times 4$ matrix with the entries of $\bs{p}^{(i)}$ on the diagonal and $0$ elsewhere.
\end{theorem}

\noindent\textsc{Proof :}  For all $n\geq 2$, and $q\geq 1$, set
\begin{equation*}
\bs{M}^{(n)}_q=\frac{1}{\sqrt{|\dT_{n-1}^*|}}\sum_{k=1}^q\left(
\begin{array}{l}
\delta_{2k}\big((1-\delta_{4k})(1-\delta_{4k+1})-p^{(0)}(0,0)\big)\\
\delta_{2k}\big(\delta_{4k}(1-\delta_{4k+1})-p^{(0)}(1,0)\big)\\
\delta_{2k}\big((1-\delta_{4k})\delta_{4k+1}-p^{(0)}(0,1)\big)\\
\delta_{2k}\big(\delta_{4k}\delta_{4k+1}-p^{(0)}(1,1)\big)\\
\delta_{2k+1}\big((1-\delta_{4k+2})(1-\delta_{4k+3})-p^{(1)}(0,0)\big)\\
\delta_{2k+1}\big(\delta_{4k+2}(1-\delta_{4k+3})-p^{(1)}(1,0)\big)\\
\delta_{2k+1}\big((1-\delta_{4k+2})\delta_{4k+3}-p^{(1)}(0,1)\big)\\
\delta_{2k+1}\big(\delta_{4k+2}\delta_{4k+3}-p^{(1)}(1,1)\big)
\end{array}\right).
\end{equation*}
Let $(\cG_q)$ be the filtration of cousin cells: $\cG_0=\sigma\{\delta_1, \delta_2, \delta_3\}$ and for all $q\geq 1$, $\cG_q=\cG_{q-1}\vee\sigma\{\delta_{4q}, \delta_{4q+1}, \delta_{4q+2}, \delta_{4q+3}\}$. For all $n\geq 2$, $(\bs{M}^{(n)}_q)$ is a $(\cG_q)$-martingale with finite moments of all order. We apply Theorem~3.II.10 of \citep{Duflo97} to this sequence of martingales and with the stopping times $\nu_n=|\dT_{n-2}|$. The $\mathbb{P}_{\overline{\cE}}$ a.s. limit of the increasing process is
\begin{equation*}
<\bs{M}^{(n)}>_{\nu_n}=  \frac{1}{|\dT^*_{n-1}|}  \left(
\begin{array}{cc}
\sum_{k \in \dT_{n-2}}\delta_{2k}\bs{V}^0&0\\
0&\sum_{k \in \dT_{n-2}}\delta_{2k+1}\bs{V}^1
\end{array}\right)
\xrightarrow[n\rightarrow\infty]{} \bs{G}=\left(
\begin{array}{cc}
z^0\bs{V}^0&0\\
0&z^1\bs{V}^1
\end{array}\right).
\end{equation*}
 In addition, the Lindeberg condition holds as the $\delta_k$ have finite moments of all order. Thus, we obtain the convergence $M^{(n)}_{|\dT_{n-2}|}
\liml
\cN(0, \bs{G} )$ on $(\overline{\cE},\mathbb{P}_{\overline{\cE}})$.
On the other hand, $\bs{\Delta}_{n-1}^{-1}|\dT^*_{n-1}|M^{(n)}_{|\dT_{n-2}|}= \sqrt{|\dT_{n-1}^*|}(\widehat{\bs{p}}_{n}-\bs{p})$, with
\begin{equation*}
\bs{\Delta}_n=\left(
\begin{array}{cc}
\sum_{\ell=1}^nZ^0_{\ell}\bs{I}_4&0\\
0&\sum_{\ell=1}^nZ^1_{\ell}\bs{I}_4
\end{array}\right),
\end{equation*}
and $\bs{I}_4$ is the identity matrix of size $4$. As $|\dT_n^*|^{-1}\sum_{\ell=1}^nZ^i_{\ell}$ converges almost surely to $z^i$ on $(\overline{\cE},\mathbb{P}_{\overline{\cE}})$, we have the asymptotic normality as announced, using Slutsky's lemma.
\hspace{\stretch{1}}$ \Box$\\

Using the asymptotic normality of the $\wh{p}_n^{(i)}(j_0,j_1)$, we can derive Wald's test for the asymmetry of the means of the reproduction laws. Set $m=\big(p^{(0)}(1,0)+p^{(0)}(0,1)+2p^{(0)}(1,1)\big)-\big(p^{(1)}(1,0)+p^{(1)}(0,1)+2p^{(1)}(1,1)\big)$ the difference of the means of the types 0 and 1 reproduction laws and $\wh{m}_n=\big(\wh{p}_n^{(0)}(1,0)+\wh{p}_n^{(0)}(0,1)+2\wh{p}_n^{(0)}(1,1)\big)-\big(\wh{p}_n^{(1)}(1,0)+\wh{p}_n^{(1)}(0,1)+2\wh{p}_n^{(1)}(1,1)\big)$ its empirical estimator. Set $\bs{H_0^{GW}}$: $m=0$ the symmetry hypothesis and $\bs{H_1^{GW}}$: $m\neq0$. Let $(Y^{GW}_n)^2$ be the test statistic defined by
$Y_n^{GW}=|\dT_{n-1}^*|^{1/2}(\wh{\Delta}_n^{GW})^{-1/2}\wh{m}_n$,
where
$\wh{\Delta}_n^{GW} = \bs{dgw}^t\wh{\bs{V}}_n\bs{dgw}$, $\bs{dgw}=(0,1,1,2,0 -1,-1,-2)^t$,
and $\wh{\bs{V}}_n$ is the empirical version of $\bs{V}$, where $z^i$ is replaced by $|\dT_n^*|^{-1}\sum_{l=1}^nZ_l^i$ and the $p^{(i)}(j_0,j_1)$ are replaced by $\wh{p}_n^{(i)}(j_0,j_1)$. Thanks to Lemma \ref {lem:CVp}, $\wh{\bs{V}}_n$ converges a.s. to $\bs{V}$ and the test statistic has the following asymptotic properties.
\begin{theorem}\label{th test GW}
Under assumption \emph{(\textbf{AO})} and the null hypothesis  $\bs{H_0^{GW}}$, one has
$(Y_n^{GW})^2
\liml
\chi^2(1)$ on $(\overline{\cE},\mathbb{P}_{\overline{\cE}})$;
and under the alternative hypothesis $\bs{H_1^{GW}}$, one has
$\lim_{n\rightarrow\infty}(Y_n^{GW})^2 = +\infty$ a.s. on $(\overline{\cE},\mathbb{P}_{\overline{\cE}})$.
\end{theorem}

\noindent\textsc{Proof :}  Let $g$ be the function defined from $\dR^8$ onto $\dR$ by 
$g(x_1,x_2,x_3,x_4,x_5,x_6,x_7,x_8)=\big(x_2+x_3+2x_4\big)-\big(x_6+x_7+2x_8\big)$,
so that $\wh{m}_n-m=g(\wh{\bs{p}}_n)-g(\bs{p})$, and $\bs{dgw}$ is the gradient of $g$. Thus, Theorem~\ref{th:TCLGW} yields
$\sqrt{|\dT_{n-1}^*|}\big(g(\wh{\bs{p}}_n)-g(\bs{p})\big)
\liml
\cN(0,\Delta^{GW})$
on $(\overline{\cE},\mathbb{P}_{\overline{\cE}})$,
with ${\Delta}^{GW} = \bs{dgw}^t{\bs{V}}\bs{dgw}$. Under the null hypothesis $\bs{H_0^{GW}}$, $g(\bs{p})=m=0$, so that
$|\dT_{n-1}^*|(\Delta^{GW})^{-1}g(\wh{\bs{p}}_n)^2
\liml
\chi^2(1)$
on $(\overline{\cE},\mathbb{P}_{\overline{\cE}})$.
One then uses Slutsky's lemma to replace $\Delta^{GW}$ by its estimator $\wh{\Delta}_n^{GW}$ and obtain the first convergence. Under the alternative hypothesis $\bs{H_1^{GW}}$, {since $Y_n^{GW}=|\dT_{n-1}^*|^{1/2}(\wh{\Delta}_n^{GW})^{-1/2}\wh{m}_n$ and $\wh{m}_n$ converges to $m \neq 0$, $Y^{GW}_n$ tends to infinity a.s. on $(\overline{\cE},\mathbb{P}_{\overline{\cE}})$} 
\hspace{\stretch{1}}$ \Box$
\section{Asymmetry in cells characteristic}
\label{s:BAR}
We now turn to the study of asymmetry of the BAR model with missing data. We first recall some asymptotic results on the estimation of the BAR parameters proved in \citep{SGM11}.

\subsection{Model and assumptions}
\label{s:BARmodel}
 We consider the first-order asymmetric BAR process given by Eq.~\reff{defbar}.
We assume that $\mathbb{E}[X_1^8]<\infty$. Moreover the parameters
satisfy
$ 0<\max(|b|, |d|) < 1.
$
Denote by $\dF=(\cF_n)$ the natural filtration associated with the BAR process: $\cF_n = \sigma\{X_k, k\in \dT_n\}$. In all the sequel, we shall make use of the following moment and independence hypotheses.
\begin{description}
\item[(AN.1)]  For all $n\geq 0$ and for all $k\in \dG_{n+1}$,
$\veps_k$ belongs to $L^8$ with
$\sup_{n\geq 0}\sup_{k\in \dG_{n+1}}
\dE[\veps_k^8|\cF_n]<\infty$ a.s.
Moreover, there exist $\sigma^2>0$ and $\rho$ such that $\forall n\geq 0$ and $k\in \dG_{n+1}$,
$\dE[\veps_k|\cF_n] = 0$, $\dE[\veps_k^2|\cF_n]=\sigma^2$ a.s.; $\forall n\geq 0$, $\forall k\neq l \in \dG_{n+1}$ with $[k/2]=[l/2]$,
$\dE[\veps_k\veps_l|\cF_n] = \rho$ {a.s.}, {where $[x]$ denotes the largest integer less than or equal to $x$.}
\item[(AN.2)] For all $n\geq 0$ the random vectors $\{(\veps_{2k} , \veps_{2k+1}), k \in\dG_n\}$ are conditionally independent given $\cF_n$.
 \item[(AI)] The sequence $(\delta_k)$ is independent from the sequences $(X_k)$ and $(\veps_k)$.
\end{description}
The least-squares estimator of $\bs{\theta}=(a,b,c,d)^t$ is given for all $n\geq 1$ by $\wh{\bs{\theta}}_n=( \wh{a}_n, \wh{b}_n, \wh{c}_n, \wh{d}_n)^t$ with
\begin{equation*}\label{defLS}
\wh{\bs{\theta}}_n = \bs{\Sigma}_{n-1}^{-1}\sum_{k \in \mathbb{T}_{n-1}}\left(
\begin{array}{c}
\delta_{2k}X_{2k}  \\
\delta_{2k}X_kX_{2k} \\
\delta_{2k+1}X_{2k+1} \\
\delta_{2k+1}X_kX_{2k+1}
\end{array}\right),\ \bs{\Sigma}_{n} = \left( \begin{array}{cc}
\bs{S}^0_{n} & 0 \\
0 & \bs{S}^1_{n}
\end{array} \right),\
\bs{S}^i_{n} =\sum_{k \in \mathbb{T}_{n}}\delta_{2k+i}\left(
\begin{array}{cc}
1 & X_k \\
X_k &X^2_k
\end{array}\right),\ \bs{S}^{0,1}_{n} =\sum_{k \in \mathbb{T}_{n}}\delta_{2k}\delta_{2k+1}\left(
\begin{array}{cc}
1 & X_k \\
X_k &X^2_k
\end{array}\right).
\end{equation*}
Denote also $\bs{L}^{0}$,  $\bs{L}^{1}$,  $\bs{L}^{0,1}$ the a.s. limits: 
$\lim_{n\rightarrow \infty} \ind{|\dG_n^*|>0}{\bs{S}^i_{n}}/{|\dT_n^*|} = \indnex\bs{L}^i$, $\lim_{n\rightarrow \infty} \ind{|\dG_n^*|>0}{\bs{S}^{0,1}_{n}}/{|\dT_n^*|} = \indnex\bs{L}^{0,1}$ (see Proposition 4.2 of \citep{SGM11}).
We now recall Theorems~3.2 and 3.4 of \citep{SGM11}.
\begin{theorem}\label{th:TCLBAR}
Under \emph{(\textbf{AN.1-2})}, \emph{(\textbf{AO})} and \emph{(\textbf{AI})}, the estimator $\widehat{\bs{\theta}}_{n}$ is strongly consistent
$\lim_{n\rightarrow\infty}\ind{|\dG_n^*|>0} \widehat{\bs{\theta}}_{n} =\ind{\overline{\cE}}\bs{\theta}$
{a.s.}
In addition, we have the asymptotic normality
$\sqrt{|\dT^*_{n-1}|} (\widehat{\bs{\theta}}_{n}-\bs{\theta})
\liml
\cN(0,\bs{\Sigma}^{-1}\bs{\Gamma} \bs{\Sigma}^{-1})$ {on }
$({\overline{\cE}}, \dP_{\overline{\cE}})$,
where 
\begin{equation*}\label{def Lambda}
\bs{\Sigma} = \left(\begin{array}{cc}\bs{L}^0&0\\0&\bs{L}^1\end{array}\right), \quad\textrm{and}\quad\bs{\Gamma} = \left(\begin{array}{cc}\sigma^2\bs{L}^0&\rho\bs{L}^{0,1}\\\rho\bs{L}^{0,1}&\sigma^2\bs{L}^1\end{array}\right).
\end{equation*}
\end{theorem}
\subsection{Asymmetry tests}
\label{s:BARtests}
Using Theorem \ref{th:TCLBAR}, we now propose two different asymmetry tests. The first one compares the couples $(a,b)$ and $(c,d)$.
Set $\bs{H_0^{c}}$: $(a,b)=(c,d)$ the symmetry hypothesis and $\bs{H_1^{c}}$: $(a,b)\neq (c,d)$. Let $(\bs{Y}_n^c)^t\bs{Y}_n^c$ be the test statistic defined by
$\bs{Y}_n^c=|\dT_{n-1}^*|^{1/2}(\wh{\bs{\Delta}}_n^c)^{-1/2}(\wh{a}_n-\wh{c}_n, \wh{b}_n-\wh{d}_n)^t$,
with
$\wh{\bs{\Delta}}_n^c = \bs{dgc}^t|\dT^*_{n-1}|{\bs{\Sigma}}_{n-1}^{-1}\wh{\bs{\Gamma}}_{n-1}|\dT^*_{n-1}|{ \bs{\Sigma}}_{n-1}^{-1}\bs{dgc}$,
\begin{equation*}
\bs{dgc}=\left(\begin{array}{cccc}1&0&-1&0\\0&1&0&-1\end{array}\right)^t,\qquad
\wh{\bs{\Gamma}}_n=\frac{1}{|\dT^*_n|}\left(\begin{array}{cc}\wh{\sigma}_{n+1}^2\bs{S}_n^0&\wh{\rho}_{n+1}\bs{S}_n^{0,1}\\\wh{\rho}_{n+1}\bs{S}_n^{0,1}&\wh{\sigma}_{n+1}^2\bs{S}_n^1\end{array}\right),
\end{equation*}
$\wh{\sigma}^2_n = { |\dT_{n}^*|^{-1}}\sum_{k \in
\dT^*_{n-1}}(\wh{\veps}_{2k}^2 +   \wh{\varepsilon}_{2k+1}^2)$, 
$\wh{\rho}_n ={|\dT_{n-1}^{*01}|^{-1}}\sum_{k \in\dT_{n-1}}
\wh{\veps}_{2k} \wh{\veps}_{2k+1}$, $\dT_{n}^{*01} =\{k \in \dT_n : \delta_{2k}  \delta_{2k+1}=1\}$ and for all $k\in\mathbb{G}_n$,
\begin{equation*}
\left\{
    \begin{array}{lcrclcl}
     \wh{\veps}_{2k} &=& \delta_{2k}(X_{2k} &-& \wh{a}_n &-& \wh{b}_{n}X_k), \vspace{1ex}\\
     \wh{\veps}_{2k+1} &=& \delta_{2k+1}(X_{2k+1} &-& \wh{c}_{n} &-& \wh{d}_{n}X_k).
    \end{array}\right.
\end{equation*}
\begin{theorem}
\label{theobar1}
Under assumptions \emph{(\textbf{AN.1-2})}, \emph{(\textbf{AO})}, \emph{(\textbf{AI})} and the null hypothesis  $\bs{H_0^{c}}$, one has
$(\bs{Y}_n^c)^t\bs{Y}_n^c
\liml
\chi^2(2)$ {on} $(\overline{\cE},\mathbb{P}_{\overline{\cE}})$;
and under the alternative hypothesis $\bs{H_1^{c}}$, one has
$\lim_{n\rightarrow\infty}\|\bs{Y}_n^c\|^2 = +\infty$ {a.s. on} $(\overline{\cE},\mathbb{P}_{\overline{\cE}})$.
\end{theorem}

\noindent\textsc{Proof :}  We mimic the proof of Theorem  \ref{th test GW} with $g$  the function defined from from $\dR^4$ onto $\dR^2$ by $g(x_1,x_2,x_3,x_4)=\big(x_1-x_3, x_2-x_4\big)^t$,
so that $\bs{dgc}$ is the gradient of $g$.
\hspace{\stretch{1}}$ \Box$\\

Our second test compares the fixed points $a/(1-b)$ and $c/(1-d)${, which are the asymptotic means of $X_{2k}$ and $X_{2k+1}$ respectively}.
Set $\bs{H_0^{f}}$: $a/(1-b)=c/(1-d)$ the symmetry hypothesis and $\bs{H_1^{f}}$: $a/(1-b)\neq c/(1-d)$. Let $(Y_n^f)^2$ be the test statistic defined by
$Y_n^f=|\dT_{n-1}^*|^{1/2}(\wh{\bs{\Delta}}_n^f)^{-1/2}\big(\wh{a}_n/(1-\wh{b}_n)-\wh{c}_n/(1-\wh{d}_n)\big)$,
where
$\wh{\bs{\Delta}}_n^f= \bs{dgf}^t|\dT^*_{n-1}|{\bs{\Sigma}}_{n-1}^{-1}\wh{\bs{\Gamma}}_{n-1}|\dT^*_{n-1}|{ \bs{\Sigma}}_{n-1}^{-1}\bs{dgf}$,
$\bs{dgf}=\big(1/(1-b), a/(1-b)^2, -1/(1-d), -c/(1-d)^2\big)^t$.
This test statistic has the following asymptotic properties.
\begin{theorem}
\label{theobar2}
Under assumptions \emph{(\textbf{AN.1-2})}, \emph{(\textbf{AO})}, \emph{(\textbf{AI})} and the null hypothesis  $\bs{H_0^{f}}$, one has
$(Y_n^f)^2
\liml
\chi^2(1)$ {on} $(\overline{\cE},\mathbb{P}_{\overline{\cE}})$;
and under the alternative hypothesis $\bs{H_1^{f}}$, one has
$\lim_{n\rightarrow\infty}(Y_n^f)^2 = +\infty$ {a.s. on} $(\overline{\cE},\mathbb{P}_{\overline{\cE}})$.
\end{theorem}

\noindent\textsc{Proof :} We mimic the proof of Theorem  \ref{th test GW} with $g$  the function defined from $\dR^4$ onto $\dR$ by 
$g(x_1,x_2,x_3,x_4)=\big(x_1/(1-x_2)-x_3/(1-x_4)\big)^t$,
so that $\bs{dgf}$ is the gradient of $g$. 
\hspace{\stretch{1}}$ \Box$

\section{Application to simulated data}
\label{s:simu}
We now study the behavior of our three tests on simulated data. For each test, we compute, in function of the generation $n$ and for different thresholds, the proportion of rejections under hypotheses $\bs{H_0}$ and $\bs{H_1}$, the latter proportion being an indicator of the power of the test. Proportions are computed on a sample of 1000 replicated trees.
\subsection{Asymmetry test for the Galton-Watson process}
In Table \ref{tGW}, we see that  the observed proportions of p-values under the thresholds (0.05, 0.01, 0.001), are 
\begin{table}[htbp]
\begin{center}
\begin{tabular}{| c |c c c |c c c  c| }
\hline
Generation & \multicolumn{3}{|c|}{Under $\bs{H_0^{GW}}$ } &  & \multicolumn{3}{c|}{Under $\bs{H_1^{GW}}$ } \\ \hline
  & $p<0.05$ &$ p<0.01$& $ p<0.001$&  & $p<0.05$ & $p<0.01$ &$ p<0.001$ \\ 
7 & 6.4 & 1.9 & 0.3 & &27.8 & 11.8 & 03.6  \\ 
8 &  5.6 & 1.4 & 0.3 & &44.2 & 22.2  & 07.6 \\ 
9 & 5.5 & 1.1 & 0.3 & & 58.6& 38.5& 17.0   \\ 
10 & 5.7 & 1.5 & 0.2 & &79.4 & 60.8 & 35.9  \\ 
11 & 4.8 & 1.0 & 0.1 & & 93.1 & 82.0 & 64.2 \\ 
\cline{1-7}
\hline
\end{tabular}
\caption{Proportions of p-values under the $0.05$, $0.01$ and $0.001$ thresholds of the asymmetry tests for the means of the GW process (1000 replicas) $\bs{p}_0=(0.04, 0.08, 0.08, 0.8)$ (under (\textbf{H1}), $\bs{p}_1=(0.15, 0.08, 0.08, 0.69)$)}
\label{tGW}
\end{center}
\end{table}
close to the expected proportions of rejection under  $\bs{H_0^{GW}}$ suggesting that the asymptotic law of the statistic $(\bs{Y}_n^{GW})^2$ is available by generation $8$. Under $\bs{H_1^{GW}}$, the  power of the test increases from $27.8$ (\%) for the generation 7 to $93.1$ (\%) for the generation 11 for a risk of type 1 fixed at $0.05$.
\subsection{Asymmetry tests for the BAR process}
The first asymmetry test compares the parameters $(a,c)$ and $(c,d)$. In Table \ref{tBAR1}, we see that  the observed proportions of p-values under the thresholds (0.05, 0.01, 0.001), are close to the expected proportions of rejection under  $\bs{H_0^{c}}$ suggesting that the asymptotic law of the statistic $\|\bs{Y}_n^c\|^2$ is available at generation $8$. Under $\bs{H_1^{c}}$, the  power of the test increases from $37.4$ (\%) for the generation 7 to $95.7$ (\%) for the generation 11 for a risk of type 1 fixed at $0.05$.
\begin{table}[htbp]
\begin{center}
\begin{tabular}{| c |c c c |c c c c| }
\hline
Gen & \multicolumn{3}{|c|}{Under $\bs{H_0^{c}}$ } &  & \multicolumn{3}{c|}{Under $\bs{H_1^{c}}$ } \\ 
\hline
  & $p<0.05$ &$ p<0.01$& $ p<0.001$&  & $p<0.05$ & $p<0.01$ &$ p<0.001$ \\ 
7 & 6.6 & 2.2 & 0.6 & & 37.4 & 19.7 & 08.0  \\ 
8 &  5.5 & 1.5 & 0.3 & & 53.6 & 31.0  & 14.6  \\ 
9 & 5.5 & 1.3 & 0.3 & & 71.1 & 52.3 & 30.3   \\ 
10 & 6.3 & 1.2 & 0.1 & & 86.8 & 75.5 & 56.1  \\ 
11 & 5.9 & 0.6 & 0.1 & & 95.7 & 90.8 & 81.4 \\ 
\cline{1-7}
\hline
\end{tabular}
\caption{Proportions of p-values under the $0.05$, $0.01$ and $0.001$ thresholds of the asymmetry test for the parameters of the BAR process (1000 replicas) $a = b = 0.5$ (under (\textbf{H1}), $c=0.5 ; d= 0.4$)}
\label{tBAR1}
\end{center}
\end{table}
\newline
For the asymmetry test for the fixed points, we see in Table \ref{tBAR2} that the observed proportions go away from the expected ones under $\bs{H_0^{f}}$ until the 10th generation, suggesting that the asymptotic law of the statistic is not reached before the 10th generation.  We also remark that the power is weak until the 10th generation. 
\begin{table}[htbp]
\begin{center}
\begin{tabular}{| c |c c c |c c  c  c| }
\hline
Gen & \multicolumn{3}{|c|}{Under $\bs{H_0^{f}} $} &  & \multicolumn{3}{c|}{Under $\bs{H_1^f}$ } \\ \hline
  & $p<0.05$ &$ p<0.01$& $ p<0.001$&  & $p<0.05$ & $p<0.01$ &$ p<0.001$ \\ 
7 & 2.2 & 0.7 & 0 & & 23.1 & 07.4 & 01.4  \\ 
8 &  3.3 & 0.5 & 0.1 & & 41.3 & 20.5  & 06.1  \\ 
9 & 3.8 & 0.5 & 0 & & 64.6 & 41.6 & 18.6   \\ 
10 & 4.7 & 0.8 & 0 & & 82.9 & 68.1 & 46.3  \\ 
11 & 5.5 & 0.7 & 0.1 & & 94.5 & 88.5 & 74.5 \\ 
\cline{1-7}
\hline
\end{tabular}
\caption{Proportions of p-values under the $0.05$, $0.01$ and $0.001$ thresholds of the asymmetry test for the fixed points of the BAR process (1000 replicas) $a = b = 0.5$ (under (\textbf{H1}), $c=0.5 ; d= 0.4$)}
\label{tBAR2}
\end{center}
\end{table}
\section{Application to real data: aging detection of Escherichia coli}
\label{s:Ecoli}
To study aging in the single cell organism E. coli, \citet{SMT05} filmed 94 colonies of dividing cells, determining the complete lineage and the growth rate of each cell. E. coli is a rod-shaped bacterium that reproduces by dividing in the middle. Each cell inherits an old end or \emph{pole} from its mother, and creates a new pole. Therefore, each cell has a \emph{type}: old pole or new pole cell inducing asymmetry in the cell division. \citet{SMT05} propose a statistical study of the genealogy and pair-wise comparison of sister cells  assuming independence between the pairs  of sister cells which is not verified in  the lineage. 

\smallskip

\begin{figure}[htbp]
 \centerline{\includegraphics[width=5cm]{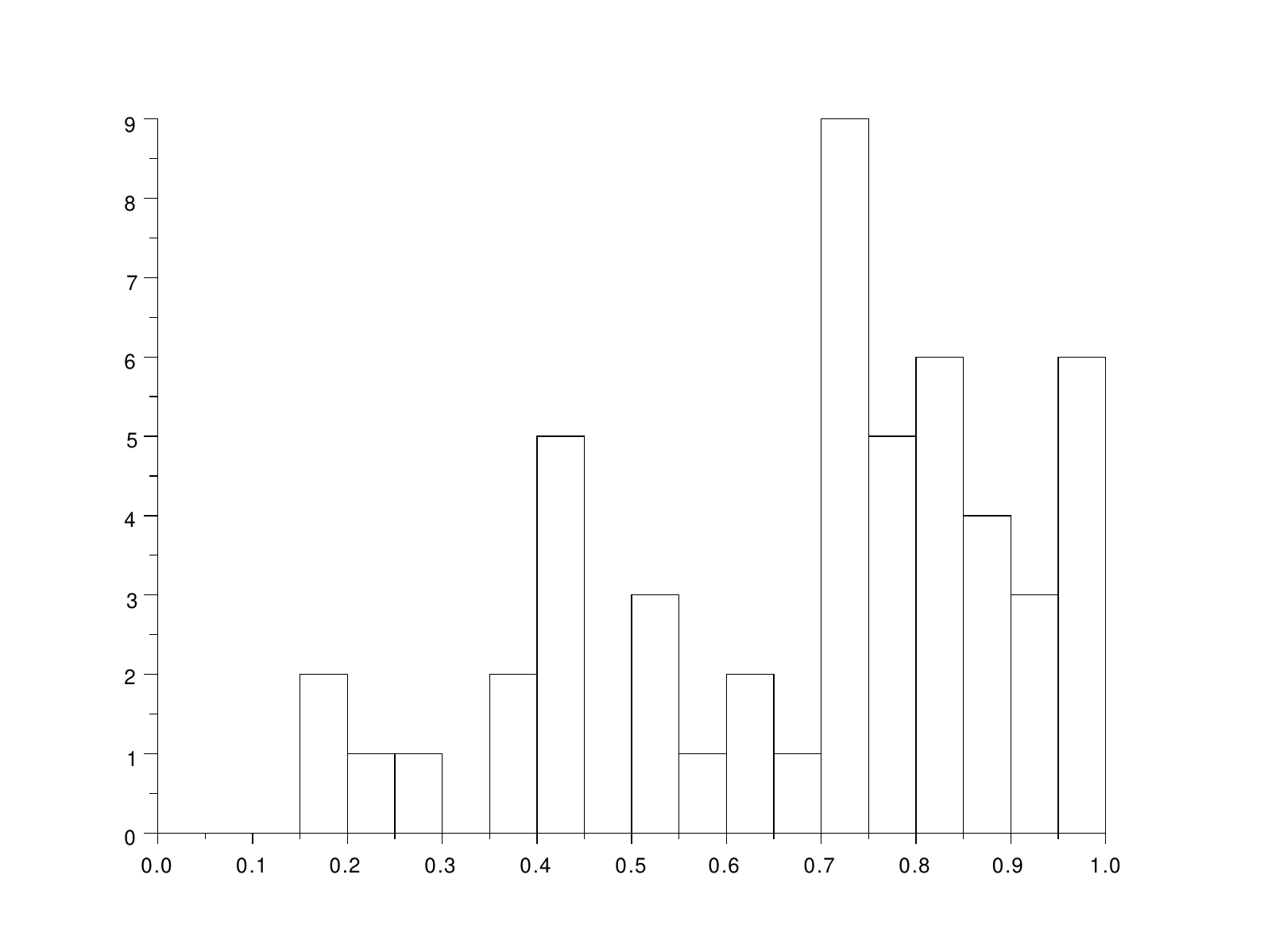}}
\caption{Histogram of the 51 p-values of the test $\bs{H_0^{GW}}$}
\label{ftGW}
\end{figure}

\begin{figure}[htbp]
  \centering
  \subfloat[$\bs{H_0^{c}}$]{\label{abcd-a}\includegraphics[width=5cm]{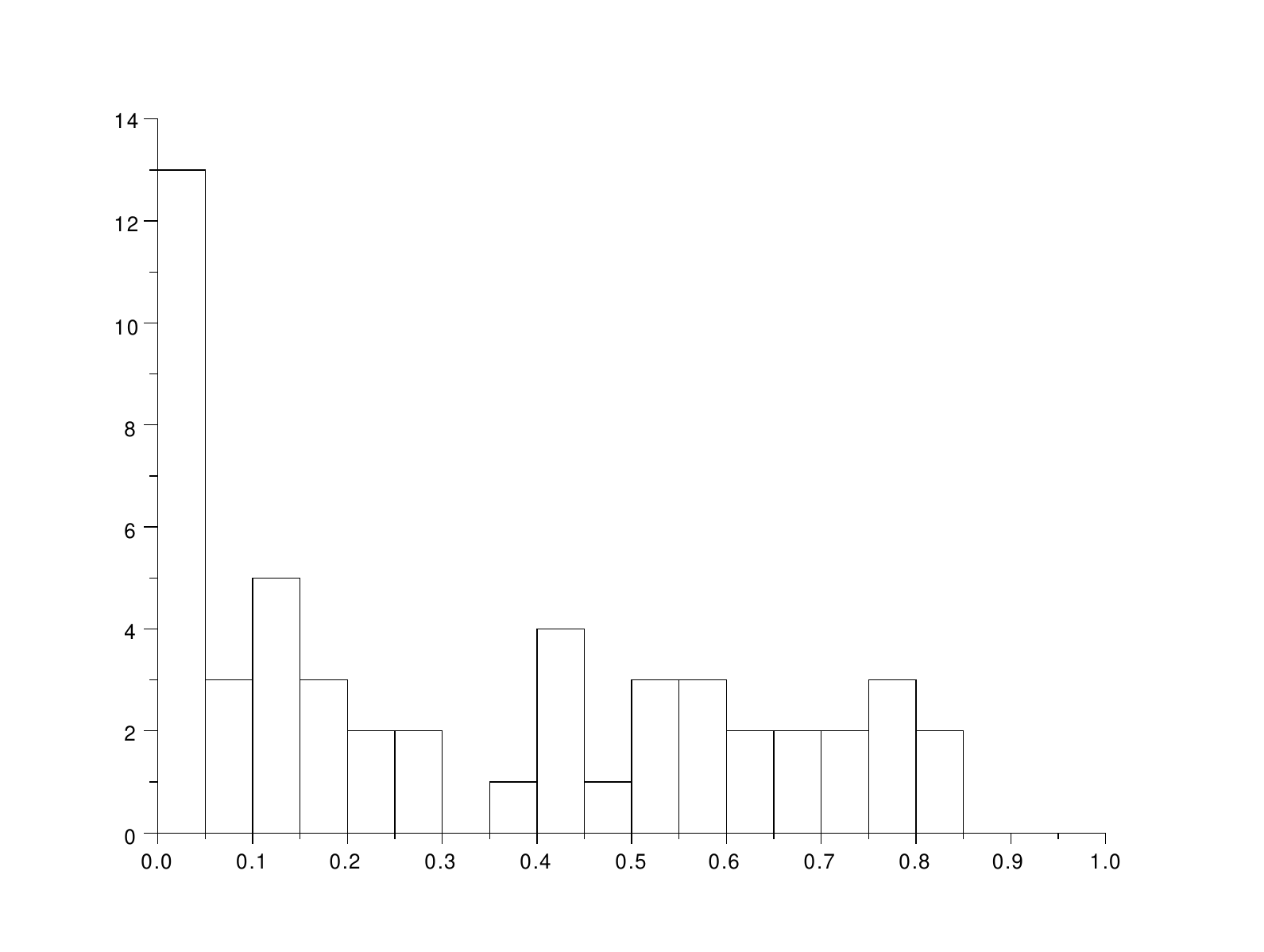}}                
  \subfloat[$\bs{H_0^{f}}$]{\label{ptfixe-b}\includegraphics[width=5cm]{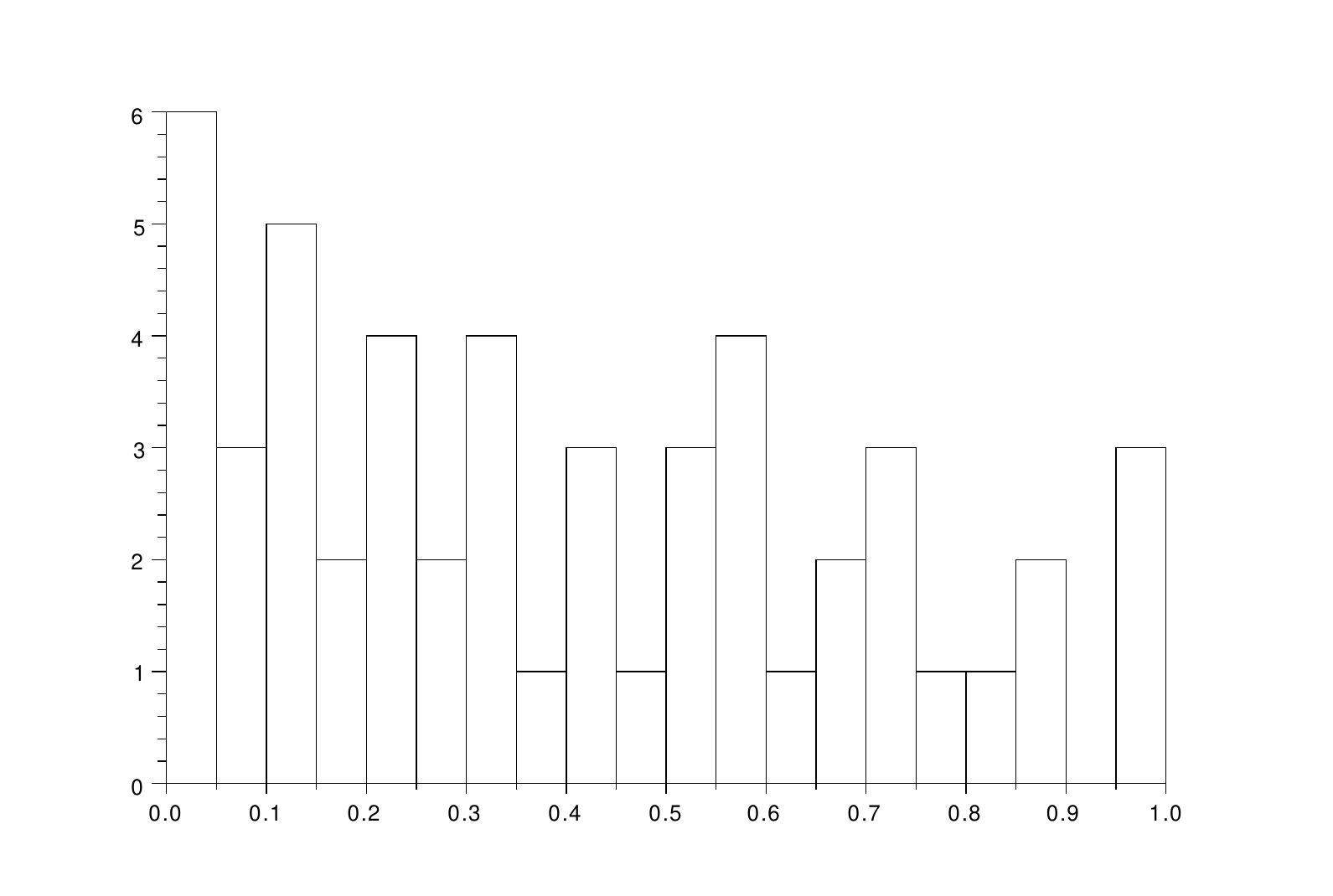}}
  \caption{Histogram of the 51 p-values of the tests of assumptions $\bs{H_0^{c}}$ and $\bs{H_0^{c}}$}
  \label{fig:tests}
\end{figure}

Figures \ref{ftGW},  \ref{abcd-a}, \ref{ptfixe-b} present the results of our tests of the null hypotheses $\bs{H_0^{GW}}$, 
$\bs{H_0^{c}}$ and $\bs{H_0^{f}}$ on the 51 data sets issued of the 94 colonies containing at least eight or nine generations. Figure \ref{ftGW} shows that the hypotheses of equality of the expected number of observed offspring between two sisters is not rejected whatever the data set. This result is not surprising: in our sets, the data are missing most frequently because the cells were out of the range of the camera.
The null hypotheses of the two tests on the BAR parameters are rejected for one set in four for $\bs{H_0^{c}}$ and for one in eight for $\bs{H_0^{f}}$.  A global conclusion on the difference between the old pole cell and the new pole cell is not easy. Regarding the results of the simulations in Tables \ref{tBAR1} and \ref{tBAR2}, this lack of evidence is probably due to a low power of the tests at generations $8$ and $9$. Some data sets with more than $9$ generations would probably show a more significant difference.

\bibliographystyle{model2-names}
 \bibliography{blabla}







\end{document}